\documentstyle{cespecial}
\input pst-plot.sty

\begin{document}
\annalsline{157}{2003}
\received{October 2, 2001}
\startingpage{689}
\def\bye{\end{document}}
 \font\tenrm=cmr10
\def\ritem#1{\item[{\rm #1}]}

\def\eqref#1{(\ref{#1})}
%--------------- Author macros ---------------
%for Bbb in amstex
\catcode`\@=11
\font\twelvemsb=msbm10 scaled 1100
\font\tenmsb=msbm10
%\font\ninemsb=msbm7 scaled 1100%msbm9
\font\ninemsb=msbm10 scaled 800
\newfam\msbfam
\textfont\msbfam=\twelvemsb  \scriptfont\msbfam=\ninemsb
  \scriptscriptfont\msbfam=\ninemsb
\def\msb@{\hexnumber@\msbfam}
\def\Bbb{\relax\ifmmode\let\next\Bbb@\else
 \def\next{\errmessage{Use \string\Bbb\space only in math
mode}}\fi\next}
\def\Bbb@#1{{\Bbb@@{#1}}}
\def\Bbb@@#1{\fam\msbfam#1}
\catcode`\@=12

 \catcode`\@=11
\font\twelveeuf=eufm10 scaled 1100
\font\teneuf=eufm10
\font\nineeuf=eufm7 scaled 1100%eufm9
\newfam\euffam
\textfont\euffam=\twelveeuf  \scriptfont\euffam=\teneuf
  \scriptscriptfont\euffam=\nineeuf
\def\euf@{\hexnumber@\euffam}
\def\frak{\relax\ifmmode\let\next\frak@\else
 \def\next{\errmessage{Use \string\frak\space only in math
mode}}\fi\next}
\def\frak@#1{{\frak@@{#1}}}
\def\frak@@#1{\fam\euffam#1}
\catcode`\@=12
%-------------- Author entries --------------------

\title{New upper bounds on sphere packings I} %Article title
\shorttitle{} % Shortened version for headline title

 \acknowledgements{Cohn was supported by an NSF Graduate Research Fellowship
and by a summer internship at Lucent Technologies, and currently
holds an American Institute of Mathematics five-year fellowship.
Elkies was supported in part by the Packard Foundation.}

 \twoauthors{Henry Cohn}{Noam Elkies}
 \institutions{Microsoft Research,
One Microsoft Way,
Redmond, WA\\
{\eightpoint {\it E-mail address\/}: cohn@microsoft.com}\\
Harvard University,
Cambridge, MA\\
{\eightpoint {\it E-mail address\/}: elkies@math.harvard.edu}
}

\newcommand{\R}{{\Bbb R}}
\newcommand{\Z}{{\Bbb Z}}
\newcommand{\Vol}{\mathop{\textup{vol}}}
\newcommand{\supp}{\mathop{\textup{supp}}}
\newcommand{\newchi}{\chi^{\phantom0}}

\centerline{\bf Abstract}
\vglue12pt
We develop an analogue for sphere packing of the linear
programming bounds for error-correcting codes, and use it to prove
upper bounds for the density of sphere packings, which are the
best bounds known at least for dimensions $4$ through~$36$. We
conjecture that our approach can be used to solve the sphere
packing problem in dimensions~8 and~24.

\def\sbni{\vglue-1pt\noindent}
\def\sni#1{\vglue-1pt\noindent{#1}. }
\def\ssni#1{\vglue-1pt\noindent\hskip18pt {#1}.}

\vglue12pt\centerline{\bf Contents}
\vglue9pt
\sni{1} Introduction
\sni{2} Lattices, Fourier transforms, and Poisson summation
\sni{3} Principal theorems
\sni{4} Homogeneous spaces
\sni{5} Conditions for a sharp bound
\sni{6} Stationary points
\sni{7} Numerical results
\sni{8} Uniqueness
\sbni Appendix A. Technicalities about density
\sbni Appendix  B. Other convex bodies
\sbni Appendix C. Numerical data
\sbni Acknowledgements
\sbni References

\vglue-12pt
\section{Introduction}
\label{sec-intro}

The sphere packing problem asks for the densest packing of
spheres into Euclidean space.  More precisely, what fraction of
$\R^n$ can be covered by congruent balls that do not intersect
except along their boundaries?  This problem fits into a broad
framework of packing problems, including error-correcting codes
and spherical codes. Linear programming bounds \cite{D} are the
most powerful known technique for producing upper bounds in such
problems.  In particular, \cite{KL} uses this technique to prove
the best bounds known for sphere packing density in high
dimensions. However, \cite{KL} does not study sphere packing
directly, but rather passes through the intermediate problem of
spherical codes. In this paper, we develop linear programming
bounds that apply directly to sphere packing, and study these
bounds numerically to prove the best bounds
known\footnote{W.-Y.~Hsiang has recently announced a solution of
the $8$-dimensional sphere packing problem \cite{Hs}, but the
details are not yet public.  His methods are apparently quite
different from ours.} for sphere packing in dimensions~$4$
through~$36$. In dimensions~$8$ and~$24$, our bounds are very
close to the densities of the known packings:  they are too high
by factors of $1.000001$ and $1.0007071$ in dimensions~$8$
and~$24$, respectively. (The best bounds previously known were
off by factors of 1.01216 and 1.27241.) We conjecture that our
techniques can be used to prove sharp bounds in~$8$ and $24$
dimensions.

The sphere packing problem in $\R^n$ is trivial for $n=1$, and
the answer has long been known for $n=2$: the standard hexagonal
packing is optimal. For $n=3$, Hales \cite{Ha} has proved that
the obvious packing, known as the ``face-centered cubic'' packing
(equivalently, the $A_3$ or $D_3$ root lattice), is  optimal, but
his proof is long and difficult, and requires extensive computer
calculation;  as of December, 2002, it has not yet been
published, but it is widely regarded as being likely to be
correct. For $n \ge 4$ the problem remains unsolved. Upper and
lower bounds on the density are known, but they differ by an
exponential factor as $n \to \infty$. Each dimension seems to
have its own peculiarities, and it does not seem likely that a
single, simple construction will give the best packing in every
dimension.

We begin with some basic background on sphere packings;  for more
information, see \cite{CS}.  Recall that a {\it lattice\/}
in~$\R^n$ is a subgroup consisting of the integer linear
combinations of a basis of~$\R^n$.  One important way to create a
sphere packing is to start with a lattice $\Lambda \subset \R^n$,
and center the spheres at the points of $\Lambda$, with radius
half the length of the shortest nonzero vectors in $\Lambda$.
Such a packing is called a {\it lattice packing\/}. Not every
sphere packing is a lattice packing, and in fact it is plausible
that in all sufficiently large dimensions, there are packings
denser than every lattice packing. However, many important
examples in low dimensions are lattice packings.

A more general notion than a lattice packing is a
{\it periodic packing\/}.
In periodic packings, the spheres are centered on the
points in the union of finitely many translates of a
lattice~$\Lambda$. In other words, the packing is still periodic
under translations by~$\Lambda$, but spheres can occur anywhere
in a fundamental parallelotope of~$\Lambda$, not just at its
corners (as in a lattice packing).

The {\it density\/} $\Delta$ of a packing is defined to be the
fraction of space covered by the balls in the packing. Density is
not necessarily well-defined for pathological packings, but in
those cases one can take a lim~sup of the densities for
increasingly large finite regions. One can prove that periodic
packings come arbitrarily close to the greatest packing density,
so when proving upper bounds it suffices to consider periodic
packings. Clearly, density is well-defined for periodic packings,
so we will not need to worry about subtleties.  See
Appendix~A for more details.

For many purposes, it is more convenient to talk about the
{\it center density\/}~$\delta$.
It is the number of sphere-centers per unit
volume, if unit spheres are used in the packing.  Thus,
$$
\Delta = \frac{\pi^{n/2}}{(n/2)!} \delta,
$$
since a unit sphere has volume $\pi^{n/2}/(n/2)!$.  Of course, for
odd $n$ we interpret $(n/2)!$  as $\Gamma(n/2+1)$.

In most dimensions, there are not even any plausible conjectures
for the densest sphere packing.  The only exceptions are low
dimensions (up to perhaps $8$ or $10$), and a handful of higher
dimensions (such as $12$, $16$, and $24$).  The most striking
examples are $8$ and $24$ dimensions. In those dimensions, the
densest packings are undoubtedly the $E_8$ root lattice and the
Leech lattice, respectively.  The $E_8$ lattice is easy to
define.  It consists of all points of $\R^8$ whose coordinates
are either all integers or all halves of odd integers, and sum to
an even integer. A more illuminating characterization is as
follows: $E_8$ is the unique lattice in $\R^8$ of covolume~$1$
such that all vectors $v$ in the lattice have even norm $\langle
v,v\rangle$. Such a lattice is called an even unimodular
lattice.  Even unimodular lattices exist only in dimensions that are multiples
of~$8$, and in $\R^8$ there is only one, up to isometries of
$\R^8$. The Leech lattice is harder to write down explicitly; see
\cite{CS} for a detailed treatment. It is the unique even
unimodular lattice in $\R^{24}$ with no vectors of length
$\sqrt{2}$. These two lattices have many remarkable properties
and connections with other branches of mathematics, but so far
these properties have not led to a proof that they are optimal
sphere packings.  We conjecture that our linear programming
bounds can be used to prove optimality.

If linear programming bounds can indeed be used to prove the
optimality of these lattices, it would not come as a complete
surprise, because other packing problems in these dimensions can
be solved similarly.  The most famous example is the kissing
problem:  how many nonoverlapping unit balls can be arranged
tangent to a given one?  If we regard the points of tangency as a
spherical code, the question becomes how many points can be
placed on a sphere with no angles less than $\pi/3$. Odlyzko and
Sloane \cite{OS} and Levenshtein \cite{Lev} independently used
linear programming bounds to solve the kissing problem in $8$ and
$24$ dimensions. (The solutions in dimensions~$8$ and $24$ are
obtained from the minimal nonzero vectors in the $E_8$ and Leech
lattices.) Because we know\break {\it a priori} that the answer must be an
integer, any upper bound within less than~$1$ of the truth would
suffice. Remarkably, the linear programming bound gives the exact
answer, with no need to take into account its integrality.  By
contrast, in most dimensions it gives a noninteger.  The
remarkable exactness seems to occur only in dimensions $1$, $2$,
$8$, and $24$.  We observe the same numerically in our case, but
can prove it only for dimension~$1$.

Figure~1 compares our results with the best
packings known as of December, 2002 (see Tables~I.1(a)
and~I.1(b) of \cite[pp.~xix, xx]{CS}), and the best upper
bounds previously known in these dimensions (due to Rogers
\cite{Ro}). The graph was normalized for comparison with Figure~15
from \cite[p.~14]{CS}.
%\figin{fig1}{980}
\begin{figure}
\begin{center}
\psset{xunit=0.25cm,yunit=1cm,linewidth=0.25pt}

\savedata{\bestknown}[{{0.0,0.0}, {1.0,-0.7604166667},
{2.0,-1.334123602}, {3.0,-1.843723030}, {4.0,-2.166666668},
{5.0,-2.510389698}, {6.0,-2.667456935}, {7.0,-2.760416667},
{8.0,-2.666666667}, {9.0,-3.093886259}, {10.0,-3.219830909},
{11.0,-3.340337786}, {12.0,-3.254772091}, {13.0,-3.340337786},
{14.0,-3.334323519}, {15.0,-3.093886259}, {16.0,-2.666666667},
{17.0,-2.760416667}, {18.0,-2.610427540}, {19.0,-2.510389698},
{20.0,-2.093093187}, {21.0,-1.843723030}, {22.0,-1.130066874},
{23.0,-0.7604166667}, {24.0,0.0}, {25.0,-0.7604100995},
{26.0,-1.334148590}, {27.0,-1.343743433}, {28.0,-1.166666667},
{29.0,-2.010410100}, {30.0,-1.875000000}, {31.0,-1.986005898},
{32.0,-1.307257948}, {33.0,-1.941891183}, {34.0,-2.389807850},
{35.0,-2.510430502}, {36.0,-2.349635295}}]

\savedata{\rogersbound}[{{0.0,0.0}, {1.0,-.7604166667},
{2.0,-1.334123602}, {3.0,-1.769134979}, {4.0,-2.096167424},
{5.0,-2.335956015}, {6.0,-2.503786868}, {7.0,-2.610793097},
{8.0,-2.666435855}, {9.0,-2.676891236}, {10.0,-2.648475338},
{11.0,-2.585754260}, {12.0,-2.493321378}, {13.0,-2.373962758},
{14.0,-2.231326547}, {15.0,-2.067841497}, {16.0,-1.885869452},
{17.0,-1.687281962}, {18.0,-1.474199786}, {19.0,-1.248144400},
{20.0,-1.010658447}, {21.0,-.7630642090}, {22.0,-.5065777090},
{23.0,-.2423650248}, {24.0,0.02854086378}, {25.0,0.3051902571},
{26.0,0.5867191273}, {27.0,0.8723058570}, {28.0,1.161190248},
{29.0,1.452654185}, {30.0,1.746043778}, {31.0,2.040729763},
{32.0,2.336135750}, {33.0,2.631697790}, {34.0,2.926904433},
{35.0,3.221263703}, {36.0,3.514299161}}]

\savedata{\newbound}[{{0.0,0.0}, {1.0,-0.7604166667},
{2.0,-1.334123602}, {3.0,-1.780494229}, {4.0,-2.096057517},
{5.0,-2.334221486}, {6.0,-2.498798538}, {7.0,-2.600839146},
{8.0,-2.649229301}, {9.0,-2.650961526}, {10.0,-2.611905963},
{11.0,-2.536974374}, {12.0,-2.430380315}, {13.0,-2.295695009},
{14.0,-2.136240609}, {15.0,-1.954425211}, {16.0,-1.752990229},
{17.0,-1.534006786}, {18.0,-1.299377863}, {19.0,-1.050913565},
{20.0,-0.7901984587}, {21.0,-0.5187296720}, {22.0,-0.2378099508},
{23.0,0.0513338332}, {24.0,0.3475636153}, {25.0,0.6501625089},
{26.0,0.9575031473}, {27.0,1.269016845}, {28.0,1.583832651},
{29.0,1.901280737}, {30.0,2.220418123}, {31.0,2.540741990},
{32.0,2.861498290}, {33.0,3.182132384}, {34.0,3.502071534},
{35.0,3.820637174}, {36.0,4.137530552}}]

\begin{pspicture}(0.0,-3.5)(36.0,4.5)
\dataplot{\bestknown} \dataplot{\rogersbound} \dataplot{\newbound}
\psaxes[Dx=4](0.0,0.0)(0.0,-3.5)(35.0,4.5)
\psaxes[labels=none,ticksize=1.5pt](0.0,0.0)(0.0,-3.5)(36.0,4.5)
\psaxes[labels=none](36.0,0.0)(36.0,-3.5)(36.0,4.5)
\rput[lB](4.0,3.5){upper curve:} \rput[lB](12.25,3.5){ Rogers'
upper bound} \rput[lB](4.0,3.0){lower curve:}
\rput[lB](12.25,3.0){ New upper bound} \rput[lB](4.0,2.5){bottom
line:} \rput[lB](12.25,2.5){ Best packing known}
% Note: the placement of the rputs (in particular, the 12.25)
% may need to be changed if the font (or font size) change.
\end{pspicture}
\end{center}
%\caption{Plot of $\log_2\delta + n(24-n)/96$ vs.\ dimension $n$.}
\centerline{Figure 1. Plot of $\log_2\delta + n(24-n)/96$ vs.\ dimension $n$.}

\label{fig-bounds}
\end{figure}

\section{Lattices, Fourier transforms, and Poisson summation}

Given a lattice $\Lambda \subset \R^n$, the
{\it dual lattice\/}
$\Lambda^*$ is defined by
$$
\Lambda^* = \{
y \ |\
\langle x,y \rangle \in \Z \textup{ for all $x \in \Lambda$}
\};
$$
it is easily seen to be the lattice with basis given by the dual
basis to any basis of $\Lambda$. The {\it covolume\/}
$|\Lambda|=\Vol(\R^n/\Lambda)$
of a lattice $\Lambda$ is the volume of any fundamental
parallelotope.  It satisfies $|\Lambda| |\Lambda^*|
 = 1$. Given any lattice $\Lambda$ with
shortest nonzero vectors of length $r$, the density of the
corresponding lattice packing is
$$
\frac{\pi^{n/2}}{(n/2)!}
\left(\frac{r}{2}\right)^n\frac{1}{|\Lambda|},
$$
and the center density is therefore $(r/2)^n/|\Lambda|$.

The Fourier transform of an $L^1$ function $f : \R^n \rightarrow
\R$ will be defined by
$$
\widehat{f}(t) = \int_{\R^n} f(x) e^{2\pi i \langle x,t \rangle}\,
dx.
$$

\proclaim{Proposition}
\label{radialft}
Let $\alpha = n/2-1$.
If $f : \R^n \rightarrow \R$ is a radial
function{\rm ,} then
$$
\widehat{f}(t) = 2\pi |t|^{-\alpha} \int_0^\infty f(r)
J_{\alpha}(2\pi r |t|) r^{n/2} \, dr,
$$
where {\rm ``}$f(r)${\rm ''} denotes the common value of $f$ on vectors of
length $r$.
\endproclaim

For a proof, see Theorem~9.10.3 of \cite{AAR}.  Here $J_{\alpha}$
denotes the Bessel function of order $\alpha$.

We will deal with functions $f : \R^n \rightarrow \R$ to which
the Poisson summation formula applies; i.e., for every
lattice $\Lambda \subset \R^n$ and every vector $v \in \R^n$,
\begin{equation}
\label{psf}
\sum_{x \in \Lambda} f(x+v) = \frac{1}{|\Lambda|}
\sum_{t \in \Lambda^*} e^{-2\pi i \langle v,t \rangle}
\widehat{f}(t),
\end{equation}
with both sides converging absolutely. It is not hard to verify
that the right-hand side of the Poisson summation formula is the
Fourier series for the left-hand side (which is periodic under
translations by elements of $\Lambda$), but of course even when
the sum on the left-hand side converges, some conditions are
needed to make it equal its Fourier series.

For our purposes, we need only the following sufficient condition:

\numbereddemo{Definition}
A function $f: \R^n \to \R$ is {\it admissible\/}
if there is a
constant $\delta>0$ such that $|f(x)|$ and $|\widehat{f}(x)|$
are bounded above by a constant times $(1+|x|)^{-n-\delta}$.
\enddemo

Admissibility implies that $f$ and $\widehat{f}$ are continuous,
and that both sides of \eqref{psf} converge absolutely.  These
two conditions alone do not suffice for Poisson summation to
hold, but admissibility does. For a proof for the integer lattice
$\Z^n$, see Corollary~2.6 of Chapter~VII of~\cite{SW}.  The
general case can be proved similarly, or derived by a linear
change of variables.

We could define admissibility more broadly, to include every
function to which Poisson summation applies, but the restricted
definition above appears to cover all the useful cases, \pagebreak and is
more concrete.

\section{Principal theorems}

Our principal result is the following theorem.  It is similar in
spirit to work of Siegel \cite{S}, but is capable of giving much
better bounds.  Gorbachev \cite{Go} has independently discovered
essentially  the same result, with a slightly different proof.  (He
concentrates on deriving Levenshtein's bound using functions $f$\break
for which $\widehat{f}$ has fairly small support, but mentions
that one could let the size of the support go to infinity.)

\proclaim{Theorem}
\label{main}
Suppose $f : \R^n \rightarrow \R$ is an admissible
function{\rm ,} is not identically zero{\rm ,} and satisfies the following
two conditions\/{\rm :}
\begin{itemize}
\ritem{(1)} $f(x)\le 0$ for $|x| \ge 1${\rm ,} and

\ritem{(2)} ${\widehat f}(t)\ge 0$ for all $t$.
\end{itemize}
Then the center density of $n$\/{\rm -}\/dimensional
sphere packings is bounded above by
$$
\frac{f(0)}{2^n{\widehat f}(0)}.
$$
\endproclaim

Notice that because $\widehat{f}$ is nonnegative and not
identically zero, we have $f(0)>0$.  If $\widehat{f}(0)=0$, then
we treat $f(0)/\widehat{f}(0)$ as $+\infty$, so the theorem is
still true, although only vacuously.

\demo{Proof}
It is enough to prove this for periodic packings, since they come
arbitrarily close to the greatest packing density (see
Appendix~A). In particular, suppose we have a
packing given by the translates of a lattice $\Lambda$ by vectors
$v_1,\dots,v_N$, whose differences are not in $\Lambda$.
If we choose the scale so that the radius of the spheres in our
packing is $1/2$ (i.e., no two centers are closer than $1$~unit),
then the center density is given by
$$
\delta = \frac{N}{2^n|\Lambda|}.
$$

By the Poisson summation formula \eqref{psf},
$$
\sum_{x \in \Lambda} f(x+v) = \frac{1}{|\Lambda|} \sum_{t \in
\Lambda^*} e^{-2\pi i\langle v,t \rangle} \widehat{f}(t)
$$
for all $v \in \R^n$.  It follows that
$$
\sum_{1 \le j,k \le N} \sum_{x \in \Lambda} f(x+v_j-v_k) =
\frac{1}{|\Lambda|} \sum_{t \in \Lambda^*}
\widehat{f}(t) \left| \sum_{1 \le j \le N} e^{2\pi i \langle v_j,
t \rangle} \right|^2.
$$
Every term on the right is nonnegative, so the sum is bounded
{}from below by the summand with $t=0$, which equals $N^2
\widehat{f}(0)/|\Lambda|$.  On the left, the vector
$x+v_j-v_k$ is the difference between two centers in the packing,
so $|x+v_j-v_k| < 1$ if and only if $x=0$ and $j=k$.  Whenever $|x+v_j-v_k|
\ge 1,$ the corresponding term in the sum is nonpositive, so we
get an upper bound of $Nf(0)$ for the entire sum.  Thus,
$$
Nf(0) \ge \frac{N^2 \widehat{f}(0)}{|\Lambda|},
$$
i.e.,
$$
\delta \le \frac{f(0)}{2^n\widehat{f}(0)},
$$
as desired.
\enddemo
\vglue9pt
This theorem was first proved by a more complicated argument,
which is given in the companion paper \cite{C}.

The hypotheses and conclusion of Theorem~\ref{main} are invariant
under rotating the function $f$. Hence, we can assume without
loss of generality that $f$ has radial symmetry, since otherwise
we can replace $f$ with the average of its rotations. The Fourier
transform maps radial functions to radial functions, and
Proposition~\ref{radialft} gives us the corresponding
one-dimensional integral transform.

As an example of how to apply Theorem~\ref{main} in one
dimension, consider the function $(1-|x|)\chi_{[-1,1]}(x)$.  It
satisfies the hypotheses of Theorem~\ref{main} in dimension
$n=1$, because it is the convolution of $\chi_{[-1/2,1/2]}(x)$
with itself, and therefore its Fourier transform is
$$
\left(\frac{\sin \pi t}{\pi t} \right)^2.
$$
Thus, this function satisfies the hypotheses of
Theorem~\ref{main}. We get a bound of $1/2$ for the center
density in one dimension, which is a sharp bound.  This example
generalizes to higher dimensions by replacing
$\chi_{[-1/2,1/2]}(x)$ with the characteristic function of a ball
about the origin.  However, the bound obtained is only the
trivial bound (density can be no greater than~$1$), so we omit
the details.  In later sections we apply Theorem~\ref{main} to
prove nontrivial bounds.

It will be useful later to have the following alternative form of
Theorem~\ref{main}:

\proclaim{Theorem}
\label{main2}
Suppose $f : \R^n \rightarrow \R$ is an admissible
function satisfying the following three conditions\/{\rm :}
\vglue4pt
{\rm (1)}\quad $f(0) = \widehat{f}(0) > 0${\rm ,}
\vglue4pt
{\rm (2)}\quad $f(x) \le 0$ for $|x| \ge r${\rm ,} and

\vglue4pt
{\rm (3)}\quad $\widehat{f}(t) \ge 0$ for all $t$.
\vglue4pt\noindent
Then the center density of sphere packings in $\R^n$ is bounded
above \pagebreak by $(r/2)^n$.
\endproclaim

Theorem 3.2 can be obtained either from rescaling the variables
in Theorem 3.1 or from the following direct proof.  For simplicity
we deal only with the case of lattice packings, but as in the
proof of Theorem 3.1 the argument extends to all periodic packings
(and hence to all packings).

\demo{{P}roof  for lattice packings}
For lattice packings, the density bound in the theorem statement
simply amounts to the claim that every lattice of covolume~$1$
contains a nonzero vector of length at most $r$.  We will prove
this first for lattices $\Lambda$ of covolume~$1-\varepsilon$,
and then let $\varepsilon \rightarrow 0+$. For such lattices,
$$
\sum_{x \in \Lambda} f(x) = \frac{1}{1-\varepsilon} \sum_{t \in
\Lambda^*} \widehat{f}(t),
$$
by Poisson summation. If all nonzero vectors in $\Lambda$ had
length greater than $r$, then all terms except $f(0)$ on the left-hand side would be nonpositive. Because all terms on
the right-hand side are nonnegative, we would have
$$
f(0) \ge \sum_{x \in \Lambda} f(x) = \frac{1}{1-\varepsilon}
\sum_{t \in \Lambda^*} \widehat{f}(t) \ge
\frac{\widehat{f}(0)}{1-\varepsilon}.
$$
However,
$$
\frac{\widehat{f}(0)}{1-\varepsilon} =
\frac{f(0)}{1-\varepsilon} > f(0),
$$
which is a contradiction.  Thus, every lattice of covolume
strictly less than~$1$ must have a nonzero vector of length $r$
or less, and it follows that the same holds for covolume~$1$.\hfill\qed
\enddemo

It seems natural to try to prove Theorem~\ref{main2} by applying
Poisson summation directly to a lattice of covolume~$1$, but some
sort of rescaling and limiting argument seems to be needed. We
included the proof to illustrate how to do this.

Logan \cite{Lo} has studied the optimization problem from
Theorem~\ref{main2} in the one-dimensional case (for reasons
unconnected to sphere packing), but we do not know of any previous
study of the higher-dimensional cases.  Unfortunately, these cases seem much more difficult than the one-dimensional
case.

\section{Homogeneous spaces}
\label{sec-homog}

The space $\R^n$ is a {\it $2$\/{\rm -}\/point homogeneous space\/};
i.e., its
isometry group acts transitively on ordered pairs of points a
given distance apart. By studying packing problems in homogeneous
spaces, one can put Theorem~\ref{main} into a broader context, in
which it can be seen to be analogous to previously known theorems
about compact homogeneous spaces.

We start by reviewing the theory of compact homogeneous spaces.
See Chapter~9 of \cite{CS} for a more detailed treatment of
this material. Suppose $X$ is a compact $2$-point homogeneous
space.   We assume that $X$ is a connected Riemannian manifold, of positive dimension.  We can write $X$ as $G/H$,
where
$(G,H)$ is a Gelfand pair of Lie groups.  Then $L^2(X)$ is a Hilbert space
direct sum of distinct irreducible representations of $G$, say
% The {} appears in the next line to keep the superscript from
% being lined up with the hat.
$\widehat{\bigoplus}{}_{j=0}^\infty V_j$. For each $j$,
evaluation gives a map $f_j : X \rightarrow V_j^*$, because $V_j$
turns out to consist of continuous functions.  We define
$$
K_j(x,y) = \langle f_j(x), f_j(y) \rangle.
$$
This is a positive definite kernel: for every finite subset $C
\subseteq X$, we have
$$
\sum_{x,y \in C} K_j(x,y)
 = \left| \sum_{x \in C}  f_j(x) \right|^2 \geq 0.
$$
Because of $G$-invariance, $K_j(x,y)$ depends only on the
distance between $x$ and $y$. This function of the distance is a
zonal spherical function; we can define a way of measuring
distance $t(x,y)$ and an ordering of the $V_j$'s so that
$K_j(x,y)$ is a polynomial $P_j$ of degree $j$ evaluated at
$t(x,y)$. In general, $t$ maps $X \times X$ to $[0,1]$, and
$t(x,y)=1$ if and only if $x=y$ (note that it is not a metric). For the unit
sphere in $\R^n$, we take $t(x,y) = (1+\langle x, y\rangle)/2$,
and the polynomial $P_j$ is the Jacobi polynomial
$P^{(\alpha,\beta)}_j(t)$, where $\alpha=\beta = (n-3)/2$.

Now suppose $C$ is a finite subset of $X$.  We get inequalities
on $C$ from the fact that for each $j$, the sum $\sum_{x \in C}
f_j(x)$ has nonnegative norm. We can apply these inequalities as
follows to get an upper bound for the size of $C$, in terms of
the minimal distance between points of $C$:

\proclaimtitle{Delsarte \cite{D}}
\proclaim{Theorem}
\label{compact}
Suppose
$$
f(t) = \sum_{j=0}^m a_j P_j(t)
$$
with $a_j \ge 0$ for all $j$ and $f(t) \le 0$ for $0 \le t \le
\tau$. If $t(x,y) \le \tau$ whenever $x$ and $y$ are distinct
points of $C${\rm ,} then
$$
|C| \le f(1)/a_0.
$$
\endproclaim

\demo{Proof}
Suppose $C$ satisfies  $t(x,y) \le \tau$ for all distinct $x,y
\in C$.  Then consider
$$
\sum_{x,y \in C} f(t(x,y)).
$$
This sum is bounded above by $|C|f(1)$ since $t(x,y) \le \tau$
unless $x=y$, and is bounded below by $|C|^2a_0$ since $f-a_0$ is
a positive definite kernel.  Thus, $|C| \le f(1)/a_0$.
\enddemo

Theorem~\ref{compact} is the analogue of Theorem~\ref{main} for
compact homogeneous spaces.  To see the analogy clearly, we need
to study $\R^n$ as a homogeneous space.

We can write $\R^n$ as $G/H$, where $G$ is the group of
isometries of $\R^n$ and $H = O(n)$.  Then we need to decompose
$L^2(\R^n)$ in terms of irreducible representations of $G$.  It
is no longer a direct sum, but it can be written as a direct
integral;  specifically, $L^2(\R^n) = \int_{0}^\infty \pi_r\,
dr$, where $\pi_r$ is the irreducible representation of $G$
consisting of functions whose Fourier transforms are
distributions with support on the sphere of radius $r$.

We can find the zonal spherical functions as follows.  The
representation~$\pi_r$ is generated by the functions $x \mapsto
e^{2\pi i\langle x,y\rangle}$ with $|y|=r$, so $\pi_r^*$ consists
of functions on the sphere of radius $r$.  The evaluation map
{}from $\R^n$ to $\pi_r^*$ takes a point $x \in \R^n$ to the
function $y \mapsto e^{2\pi i\langle x, y\rangle}$ on the sphere
of radius $r$.  Thus, the zonal spherical functions are given by
$$
K_r(x_1,x_2) = \int_{|y|=r} e^{2\pi i \langle y, x_1-x_2 \rangle} \, dy.
$$
(This of course depends only on $|x_1-x_2|$, and can be evaluated
explicitly in terms of Bessel functions using
Proposition~\ref{radialft}.) In other words, they are given by
functions whose Fourier transforms are delta functions on spheres
centered at the origin. See Section~4.15 of \cite{DM} for a more
detailed discussion of this point of view.

Now the analogue of positive combinations of the zonal spherical
functions $P_j(t)$ from the compact case is radial functions with
nonnegative Fourier transform, and we can see that
Theorem~\ref{main} corresponds to~\ref{compact}.

\vglue-6pt
\section{Conditions for a sharp bound}
\label{sec-exact}
\vglue-6pt

In one dimension, we have already seen how to use
Theorem~\ref{main} to solve the (admittedly trivial) sphere
packing problem.  Based on numerical evidence and analogy with
the kissing problem, we conjecture that it can also be used to
get sharp bounds in dimensions~$2$, $8$, and $24$.  For reasons
to be explained shortly, it is more convenient to work with
Theorem~\ref{main2} instead of Theorem~\ref{main}, so we shall do
so; we can convert everything to the framework of
Theorem~\ref{main} by rescaling the variables.

In each of dimensions $1$, $2$, $8$, and $24$, the densest known
packing is a lattice packing, given by a lattice that is
homothetic to its dual. This lattice is $\Z$ in dimension~1, the
$A_2$ root lattice (i.e., the hexagonal lattice) in dimension~2,
the $E_8$ root lattice in dimension~8, and the Leech lattice in
dimension~24.  See \cite{CS} for information about these
lattices.  Each of these lattices except $A_2$ actually equals
its dual, but that is not true for $A_2$.  However, we can
rescale $A_2$ so that the rescaled lattice is {\it isodual\/},
i.e., isometric with its own dual (in this case, via a rotation).

Suppose $\Lambda$ is any lattice of covolume~$1$, such as an
isodual lattice,  and $f$ is a radial function giving a sharp
bound on $\Lambda$ via Theorem~\ref{main2} (i.e., $r$ is the
length of the shortest nonzero vector of $\Lambda$). By Poisson
summation, we have
$$
\sum_{x \in \Lambda} f(x) = \sum_{x \in \Lambda^*} \widehat{f}(x).
$$
Given the inequalities on $f$ and $\widehat{f}$, the only way
this equation can hold is if $f$ vanishes on
$\Lambda\setminus\{0\}$ and $\widehat{f}$ vanishes on
$\Lambda^*\setminus\{0\}$.  This puts strong constraints on $f$
and $\widehat{f}$.  When $\Lambda$ is isodual, the vector lengths
in $\Lambda$ and $\Lambda^*$ are the same, so $f$ and
$\widehat{f}$ must both vanish on $\Lambda\setminus\{0\}$.

Of course, there are similar constraints on $f$ for a sharp bound
in Theorem~\ref{main} (as opposed to Theorem~\ref{main2}), but we
prefer to work with this context, since the isodual
normalizations are more pleasant, and are the standard
normalizations for $E_8$ and the Leech lattice.

It is natural to try to guess $f$ from our knowledge of its roots.
For example, in one dimension we could try
$$
f(x) = (1-x^2) \prod_{k \ge 2} \left(1 - \frac{x^2}{k^2}\right)^2
= \frac{1}{1-x^2}\left(\frac{\sin \pi x}{\pi x}\right)^2,
$$
which clearly satisfies $f(x) \le 0$ for $|x| \ge 1$ and has the
right zeros. In fact, one can compute its Fourier transform and
check that $\widehat{f}$ is nonnegative everywhere (it has
support $[-1,1]$ and is positive in $(-1,1)$), so it solves the
sphere packing problem in dimension $1$, in a different way from
the function in the previous section.

Unfortunately, it seems difficult to generalize this approach to
higher dimensions.  One can generalize this function by replacing
the sine function with a Bessel function (see
Proposition~\ref{levensh}), but that does not yield a sharp bound
in dimensions greater than~$1$.  Attempts to write down a product
with zeros at the right places for a sharp bound lead to products
that seem intractable.

One important thing to note is that for a sharp bound above
dimension~1, it is not possible for $\widehat{f}$ to have compact
support, as it does in the examples involving sine and Bessel
functions. If it did, then $f$ could not have sufficiently
densely-spaced zeros. To be precise, if $\widehat{f}$ is a radial
function with support in the ball $B(0,R)$ of radius $R$ about
the origin, then the common value $f(r)$ on vectors of radius $r$
satisfies
$$
f(r) = \int_{B(0,R)} \widehat{f}(t) e^{2\pi  i \langle rx, t
\rangle} \, dt,
$$
where $x$ is any vector with $|x|=1$. This defines an entire
function of $r$, and for all complex $r$,
$$
|f(r)| \le e^{2\pi R |r|} \int_{B(0,R)} |\widehat{f}(t)| \, dt,
$$
so $f$ is a function of exponential type, and Jensen's formula
implies that $f$ can have at most linearly spaced zeros (see
Section~15.20 of \cite{Ru}).  However, the nonzero vectors in
the Leech lattice have lengths $\sqrt{2k}$ for integers $k>1$,
and those in $E_8$ have lengths $\sqrt{2k}$ for integers $k>0$.
The function $f$ must vanish at those vector lengths, and these
roots are too densely spaced for $\widehat{f}$ to have compact
support.  Of course, $f$ also cannot have compact support
(because $\widehat{f}$ vanishes on $\Lambda^*\setminus\{0\}$).

One might wonder whether the restriction to radial functions is
misleading: perhaps a nonradial function could be constructed
more naturally.  We cannot rule out that possibility, but consider
it unlikely.  Even if $f$ is not radial, a sharp bound implies
that $f$ and $\widehat{f}$ must vanish on concentric spheres
centered at the origin and passing through the nonzero points of
$\Lambda$ and $\Lambda^*$, respectively.  The simplest reason is
that if $f$ proves that $\Lambda$ is optimal, then it proves the
same for every rotation of $\Lambda$.  Alternatively, after
rotational symmetrization $f$ and $\widehat{f}$ must vanish on
these spheres, and the inequalities on their values then imply
that they must have vanished before symmetrization (the average
of nonnegative values vanishes if and only if the values all do).
It would seem strange for $f$ to vanish on these spheres without
being radial, but of course we cannot rule it out.

\vglue-6pt
\section{Stationary points}
\label{sec-stationary}
\vglue-6pt

We do not know how to use Theorem~\ref{main} to match the best
density bound known in high dimensions, that of Kabatiansky and
Levenshtein \cite{KL}. However, it provides a new proof of the
second-best bound known, due to Levenshtein \cite{Lev}:
$$
\Delta \le \frac{j_{n/2}^n}{(n/2)!^2 4^n},
$$
where $j_t$ is the smallest positive zero of the Bessel function
$J_t$.  (For more information about the asymptotics of this bound
and how it compares with other bounds, see page~19 of
\cite{CS}, but note that equation~(42) is missing the exponent
in $j_{n/2}^n$.) We will show how to use a calculus of variations
argument to find functions that prove that bound.  This approach
is analogous to that used by Levenshtein.  Yudin \cite{Y} has
also given a proof of Levenshtein's bound that seems reminiscent
of our general approach, but not identical.

To construct a function $f$ for use in Theorem~\ref{main}, we
begin by supposing that there is a function $g$ such that $f(x) =
(1-|x|^2)\widehat{g}(x)^2$, so that $f$ automatically satisfies
the inequality $f(x) \le 0$ for $|x| \ge 1$. (We write
$\widehat{g}$ instead of $g$ for convenience later.) Assume that
$g$ is radial, and has support in the ball of radius $R$ about
the origin;  we discuss these assumptions later. Notice that
nothing in our setup requires $\widehat{f}$ to be nonnegative,
so we must check for this property later.

We have $\widehat{f} = g * (g+Lg)$, where $*$ denotes convolution
and
$$
L = \frac{1}{4\pi^2}\sum_{j=1}^n \frac{\partial^2}{\partial
t_j^2},
$$
so that under the Fourier transform, $L$ corresponds to
multiplication by $-|x|^2$. We require $f(0)=1$, i.e.,
$\int_{\R^n} g = 1$. We want to maximize $\widehat{f}(0)$,
subject to this constraint. Notice that
$$
\widehat{f}(0) = \int_{\R^n} g(0-t) (g+Lg)(t) \, dt
= \int_{\R^n} g(t) (g+Lg)(t) \, dt = \int_{\R^n} g(g+Lg),
$$
because $g$ is radial (and hence even).

Perturb $g$ to $g+h$, where $h$ has integral zero (so that $f(0)$
does not change) and has support in the ball of radius $R$.  Then
the first order change in $\widehat{f}(0)$ is
$$
\int_{\R^n} (gh+hg + gLh + hLg) = 2\int_{\R^n} h(g+Lg)
$$
(where the equality comes from integration by parts).

In order to have this vanish whenever $\int_{\R^n} h = 0$, the
function $g+Lg$ must be constant (within the support of $g$), so
for some constant $c$ we have
$$
g+Lg = c \newchi_R,
$$
where $\newchi_R$ is the characteristic function of the ball of
radius $R$ about $0$. Then it follows from
Proposition~\ref{radialft} that
$$
(1-|x|^2)\widehat{g}(x) = c (R/|x|)^{n/2} J_{n/2}(2\pi R |x|).
$$

Now $2\pi R$ must be a zero of $J_{n/2}$ for the right-hand side
to vanish at $|x|=1$, and $c$ is determined (given $R$) by
$$
1 = \int_{\R^n} g = \int_{\R^n} g+Lg =
\int_{\R^n} c \newchi_R = c R^n \pi^{n/2}/(n/2)! \, .
$$
(The function $g+Lg$ has the same integral as $g$, because
integration by parts shows that $Lg$ has integral $0$.)

In fact, we take $2\pi R$ to be the first positive root $j_{n/2}$
of $J_{n/2}$, in order to make $\widehat{f}$ nonnegative
everywhere.  To check that it is indeed nonnegative everywhere
then, we use the equation
$$
\widehat{f} = g * (g+Lg) = g * (c \newchi_R).
$$
If $g$ is nonnegative everywhere, then so is $\widehat{f}$.  We
can now determine $g$ explicitly:

We know that $g+Lg = c\newchi_R$.  The differential operator
$1+L$ takes a radial function $u(|t|)$ to
$$
\frac{u''(|t|)}{4\pi^2} + \frac{n-1}{4\pi^2|t|} u'(|t|) + u(|t|).
$$
It follows from this and the differential equation satisfied by
the Bessel functions that
$$
J_{\alpha}(2\pi |t|)/|t|^{\alpha}
$$
is in the kernel of $1+L$, if $\alpha = n/2-1$.
{}From that, one can check that
$$
g(t) =
\left(-\frac{R^{\alpha}}{J_{\alpha}(2\pi R)}
\frac{J_{\alpha}(2\pi|t|)}{|t|^{\alpha}}+1 \right)
\newchi_R(t).
$$
One can check that the minimum of
$J_{\alpha}(2\pi |t|)/|t|^{\alpha}$
occurs at
$$
|t|=j_{n/2}/(2\pi) = R,
$$
so $g$ is nonnegative everywhere, as desired,
since $J_{\alpha}(2\pi R) = J_{\alpha}(j_{n/2}) < 0$
(the roots of the functions $J_{\alpha}$ and
$J_{\alpha+1}=J_{n/2}$ are interlaced).

Thus, we have constructed a function $f$ satisfying the
hypotheses of Theorem~\ref{main}.  It has $f(0)=1$ and
$$
\widehat{f}(0)= \int_{\R^n} g(g+Lg) = c \int_{\R^n} g = c =
\frac{(n/2)!}{\pi^{n/2}R^n}.
$$
Because $R = j_{n/2}/(2\pi)$, this function proves Levenshtein's
sphere packing bound:

\proclaim{Proposition}
\label{levensh}
The function
$$
f(x) = \frac{J_{n/2}(j_{n/2}|x|)^2}{(1-|x|^2)|x|^n}
$$
satisfies the hypotheses of Theorem~{\rm \ref{main},} and leads to the
upper bound
$$
\frac{j_{n/2}^n}{(n/2)!^2 4^n}
$$
for the densities of $n$-dimensional sphere packings.
\endproclaim

This function does not optimize the bound in Theorem~\ref{main},
but it does optimize it within the class of functions whose
Fourier transforms have support in the ball of radius
$j_{n/2}/\pi$ about the origin.  This was first proved by
Gorbachev \cite{Go}.  For another proof, see \cite{C}.

It is not a coincidence that this proves exactly the same bound
as in Levenshtein's paper \cite{Lev}.  Levenshtein studies
spherical codes with minimal angular separation $\theta$, and
derives his sphere packing bound from letting $\theta \rightarrow
0$.  Under that limit, the functions he uses in the linear
programming bounds for spherical codes become our Bessel function
example.  This can be proved using the limit
\begin{equation}
\label{bmp}
\lim_{j \rightarrow \infty} j^{-\alpha}
P_j^{(\alpha,\beta)}\left(\cos \frac{z}{j}\right) =
(z/2)^{-\alpha} J_\alpha(z),
\end{equation}
which is 10.8~(41) in \cite{EMOT}.

The functions we have obtained are not optimal in any dimension
above~$1$.  There are two reasons for this.  First, we restricted
our attention to functions such that $\widehat{f}$ has compact
support, and as we have seen in Section~\ref{sec-exact}, that
cannot be true if we are to get sharp bounds. Second, and more
importantly, we implicitly considered only functions such that
$\widehat{f}$ is positive within its support.  The problem is
that if $\widehat{f}$ vanishes somewhere, then the perturbations
$h$ must be chosen so as not to push $\widehat{f}$ below zero. We
made no attempt to do so, and therefore could not find any
stationary points for which $\widehat{f}$ has zeros within its
support.  We have seen in Section~\ref{sec-exact} that zeros are
essential for sharp bounds.

Unfortunately, it seems difficult to adapt the stationary point
argument to deal with these difficulties.  One approach is to
assume that $\widehat{f}$ has zeros at certain locations, and
look at only the perturbations $h$ that, up to first order, do
not push $\widehat{f}$ below zero.  Although we can set up such
problems, we have not been able to solve them.

\section{Numerical results}
\label{sec-numerical}

It is possible to get numerical results by using linear
programming to find functions for use in Theorem~\ref{main}, as
was done for the kissing problem by Odlyzko and Sloane in
\cite{OS}.  The idea is to fix one of $f(0)$ and $\widehat{f}(0)$,
and view extremizing the other as an infinite-dimensional linear
programming problem.  One can try to approximate it with a
finite-dimensional problem, and solve it on a computer.  Although
we obtained some numerical results this way, it was cumbersome and
generally ineffective. Instead, we use the following approach.

First, consider trying to use our techniques to bound the density
of an isodual lattice.  There is no reason for optimal sphere
packings to be isodual lattices, and for example in three
dimensions they are known not to be, but it is convenient to use
this case as a stepping stone.

\proclaim{Proposition}
\label{selfdual}
Suppose $g : \R^n \rightarrow \R$ is a radial{\rm ,}
admissible function{\rm ,} is not identically zero{\rm ,} and satisfies the
following three properties\/{\rm :}\/
\begin{itemize}
\ritem{(1)} $g(0)=0${\rm ,}
\ritem{(2)} $g(x) \ge 0$ for $|x| \ge r${\rm ,} and
\ritem{(3)} $\widehat{g}=-g$.
\end{itemize}
Then every isodual lattice in dimension~$n$ must contain a
nonzero vector of length at most $r$.
\endproclaim

\demo{Proof  of special case}
For simplicity, we deal only with the case in which $g(x) > 0$
for $|x| \gg 0$. Let $\Lambda \cong \Lambda^*$ be an isodual
lattice. By Poisson summation,
$$
\sum_{x \in \Lambda} g(x) = \sum_{x \in \Lambda^*} \widehat{g}(x)
= -\sum_{x \in \Lambda} g(x),
$$
so
$$
\sum_{x \in \Lambda} g(x) = 0.
$$
In order for the sum not to be positive, the lattice $\Lambda$
must contain some nonzero vector of length at most $r$.
\enddemo

We can actually remove the hypothesis that $g(x) > 0$ for $|x|
\gg 0$ from the proof of Proposition~\ref{selfdual}, by using a
scaling trick, as in the proof of Theorem~\ref{main2}.  However,
we omit the details, because the hypothesis holds in all our
numerical examples.

Notice that given any function $f$ that proves a bound in
Theorem~\ref{main2}, we can produce a $g$ that proves the same
bound in Proposition~\ref{selfdual}, by taking $g =
\widehat{f}-f$. Thus, the bound we get for isodual lattices is at
least as good as that for general sphere packings.  In principle,
it could be better, but in practice we find that it is not (as we
explain below).

We can find excellent functions for use in
Proposition~\ref{selfdual} as follows.  Let $L_k^{\alpha}(x)$ be
the Laguerre polynomials orthogonal with respect to the measure
$e^{-x} x^{\alpha}\, dx$ on $[0,\infty)$.  Set $\alpha = n/2-1$,
and define
\vglue4pt \centerline{$
g_k(x) = L_k^{\alpha}(2\pi|x|^2)e^{-\pi|x|^2}
$}
\vglue4pt\noindent
for $k \ge 0$; we suppress the dependence on $n$ in our notation.
These functions form a basis for the radial eigenfunctions of the
Fourier transform, with eigenvalues $(-1)^k$ (see Section~4.23
and equation~(4.20.3) of \cite{Leb}).

To find a function $g$ for use in Proposition~\ref{selfdual}, we
consider a linear combination of $g_1,g_3,\dots,g_{4m+3}$, and
require it to have a root at $0$ and $m$ double roots at
$z_1,\dots,z_m$. (Counting degrees of freedom suggests that there
should be a unique such function, up to scaling.) We then choose
the locations of $z_1,\dots,z_m$ to minimize the value $r$ of the
last sign change of $g$.  To make this choice, we do a computer
search.  Specifically, we make an initial guess for the locations
of $z_1,\dots,z_m$, and then see whether we can perturb them to decrease~$r$.  We repeat the perturbations until
reaching a local optimum.  Strictly speaking, we cannot prove that it ever
converges, or comes anywhere near the global optimum.  However,
it works well in practice. At any rate this cannot affect the
validity of our bounds, only their optimality given~$m$. As we
increase $m$, this method should give better and better bounds,
which should converge to the best bounds obtainable using
Proposition~\ref{selfdual}.

This method gives good functions for use in
Proposition~\ref{selfdual}, but naturally bounds on all sphere
packings would be better than bounds only on isodual lattices. We
can turn these functions into functions satisfying the hypotheses
of Theorem~\ref{main2}, without changing $r$, as follows.

Let $h$ be a linear combination of $g_0, g_2,\dots, g_{4m+2}$
with double zeros at $z_1,\dots,z_m$, such that $g+h$ has a
double zero at $r$.  Then in the examples we have computed, $g+h$
has constant sign, which we can take to be positive. We end up
with a function $f= -g+h$ that is nonpositive outside radius
$r$, and whose Fourier transform $\widehat{f}=g+h$ is nonnegative
everywhere; furthermore, $f(0)=\widehat{f}(0)$.  Thus, we can
apply Theorem~\ref{main2} to get the same bound for general
sphere packings that $g$ proves for isodual lattices.  Notice
that it is {\it not} clear {\it a priori} that $f$ must satisfy all
the hypotheses of Theorem~\ref{main2}, but it happens in all of
our numerical examples.

Figure~1 (from \S\ref{sec-intro}) and
Table~3 (from Appendix~C)
illustrate the bounds this method produces, using $m=6$ as the
number of forced double zeros (except in dimension $7$ and lower,
where $m=5$ suffices for the accuracy we desire). We also list
for comparison Rogers' bound \cite{Ro}, which was previously the
best bound known in dimensions 4 through 36. The choices of
forced double roots are described in Table~4 (from
Appendix~C).

These bounds were calculated using a computer.  However, the
mathematics behind the calculations is rigorous.  In particular,
we use exact rational arithmetic, and apply Sturm's theorem to
count real roots and make sure we do not miss any sign changes.
More precisely, we take the quantities $2\pi z_i^2$ to be rational
numbers.  That is convenient, because the functions $g_k(x)$ are
polynomials in $2\pi |x|^2$ with rational coefficients (times
Gaussians, which have no sign changes). The one
subtlety is that in constructing $h$, we want $g+h$ to have a
double root at $r$, and $2\pi r^2$ is generally not rational.
Instead of doing this computation using floating point
arithmetic, we replace $2\pi r^2$ by a nearby rational number and
then determine $h$ exactly, using that value of $r$.  To do so,
we must increase $r$ slightly, but of course we can make the
increase as small as we wish.

In dimensions $8$ and $24$, we carried out the calculations for
$m=11$.  The resulting upper bounds are within factors of
$1.000001$ and $1.0007071$ of equality, respectively.
More precisely, in dimension~8 we take
$$
2\pi r^2 =
12.56637375,
$$
and in dimension~24 we take
$$
2 \pi r^2 =
25.1342216.
$$
The forced double roots we
used to achieve these bounds are given in Tables~1
and~2.
$$
\begin{tabular}{c|cccccc}
$i$ & 1 & 2 & 3 & 4 & 5 & 6 \\ \hline
$2\pi z_i^2$ &
$25.1328$ & $37.6995$ & $50.2678$ & $62.8463$ & $75.4579$
& $88.2463$
\end{tabular}
$$
$$
\begin{tabular}{c|ccccc}
$i$ & 7 & 8 & 9 & 10 & 11 \\ \hline
$2\pi z_i^2$ &
$101.3618$ & $115.6443$ & $131.0298$ & $150.0861$ & $174.2876$
\end{tabular}
$$
\centerline{Table 1. Forced double roots for $m=11$, $n=8$.}

% $$\qquad
%\begin{tabular}{c|cccccc}
%$i$ & 1 & 2 & 3 & 4 & 5 & 6 \\ \hline
%$2\pi z_i^2$ &
%$25.1328$ & $37.6995$ & $50.2678$ & $62.8463$ & $75.4579$
%& $88.2463$
%\end{tabular}
%$$
%$$
%\qquad\begin{tabular}{c|ccccc}
%$i$ & 7 & 8 & 9 & 10 & 11 \\ \hline
%$2\pi z_i^2$ &
%$101.3618$ & $115.6443$ & $131.0298$ & $150.0861$ & $174.2876$
%\end{tabular}
%$$
%\centerline{Table 1. Forced double roots for $m=11$, $n=8$.}

$$
\begin{tabular}{c|cccccc}
$i$ & 1 & 2 & 3 & 4 & 5 & 6 \\ \hline
$2\pi z_i^2$ &
$37.705$ & $50.285$ & $62.893$ & $75.578$ &
$88.454$ & $101.737$
\end{tabular}
$$
$$
\begin{tabular}{c|ccccc}
$i$ & 7 & 8 & 9 & 10 & 11 \\ \hline
$2\pi z_i^2$ &
$115.776$ & $131.035$ & $148.162$ & $168.215$ & $193.766$
\end{tabular}
$$
\centerline{Table 2. Forced double roots for $m=11$, $n=24$.}

%$$
%\qquad\begin{tabular}{c|cccccc}
%$i$ & 1 & 2 & 3 & 4 & 5 & 6 \\ \hline
%$2\pi z_i^2$ &
%$37.705$ & $50.285$ & $62.893$ & $75.578$ &
%$88.454$ & $101.737$
%\end{tabular}\qquad
%$$
%
%$$\qquad \quad\
%\begin{tabular}{c|ccccc}
%$i$ & 7 & 8 & 9 & 10 & 11 \\ \hline
%$2\pi z_i^2$ &
%$115.776$ & $131.035$ & $148.162$ & $168.215$ & $193.766$
%\end{tabular}
%\begin{tabular}{c}\kern4in\end{tabular}
%$$
%\centerline{Table 2. Forced double roots for $m=11$, $n=24$.}

\vglue12pt
Our numerical results lead us to the following conjectures.

\proclaim{{C}onjecture}
\label{equalconj}
The smallest possible value of $r$ in
Proposition~{\rm \ref{selfdual}} equals that in Theorem~{\rm \ref{main2},}
and for each optimal $g$ from Proposition~{\rm \ref{selfdual},} there
exists an optimal $f$ from Theorem~{\rm \ref{main2}} such that $g =
\widehat{f} -f $.
\endproclaim

\proclaim{{C}onjecture}
\label{strongconj}
There exist functions that satisfy the
hypotheses of Theorem~{\rm \ref{main2}} and solve the sphere packing
problem in dimensions $2${\rm ,} $8${\rm ,} and $24$.
\endproclaim

\proclaim{{C}onjecture}
\label{convconj}
The numerical method described above gives
bounds that converge {\rm (}\/as $m \rightarrow \infty${\rm )} to the optimal
bounds obtainable using Theorem~{\rm \ref{main2}.}
\endproclaim

\section{Uniqueness}

It is natural to ask whether the densest sphere packing in $\R^n$
is unique.  Of course, it is trivially not unique, since for
example removing a single sphere does not change the global
density. However, it is conjectured that $E_8$ and the Leech
lattice are unique among periodic packings.  By contrast, in
$\R^3$ there are infinitely many distinct periodic packings of
maximal density. Our techniques can be used to prove this
uniqueness in $8$ and $24$ dimensions, given the following slight
strengthening of Conjecture~\ref{strongconj}.  At the same time,
we will deal with the hexagonal lattice, although uniqueness in
that case has long been known.

Let $\Lambda_2$, $\Lambda_8$, and $\Lambda_{24}$ denote the
isodual scaling of the hexagonal lattice, the $E_8$ root lattice,
and the Leech lattice, respectively.

\proclaim{{C}onjecture}
\label{strongerconj} For $n \in \{2,8,24\}${\rm ,} there exists a
function that satisfies the hypotheses of Theorem~{\rm \ref{main2}} to
prove that $\Lambda_n$ is the densest packing in $\R^n$.
Furthermore{\rm ,} this function and its Fourier transform have roots
only at the vector lengths in $\Lambda_n$.
\endproclaim

First, let $n$ be $8$ or $24$, so that $\Lambda_n$ is an even
unimodular lattice (we will discuss the $n=2$ case below). Let
$f$ be a function satisfying the conclusions of
Conjecture~\ref{strongerconj}. Suppose we have a maximally dense
packing given by the translates of a lattice $\Lambda$ by vectors
$v_1,\dots,v_N$, whose differences are not in $\Lambda$.
\eject
\noindent Without
loss of generality, we can assume that $|\Lambda|=N$ and $v_1=0$.
Note that $|\Lambda|=N$ implies that the packing uses balls of
the same radius as those in $\Lambda_n$.  By Poisson summation,
$$
\sum_{1 \le j,k \le N} \sum_{x \in \Lambda} f(x+v_j-v_k) =
\frac{1}{|\Lambda|} \sum_{t \in \Lambda^*} \widehat{f}(t) \left|
\sum_{1 \le j \le N} e^{2\pi i \langle v_j, t \rangle} \right|^2,
$$
and the $f(0)$ terms cancel the $\widehat{f}(0)$ term, so we can
draw the usual conclusions from the inequalities on each side. In
particular, each vector $x+v_j-v_k$ must occur at a root of $f$.
Now we can apply the following lemma:

\proclaim{Lemma}
\label{evengen} Suppose $S$ is a subset of $\R^n$ such that $0
\in S${\rm ,} there are $n$ linearly independent vectors in $S${\rm ,} and
for all $x,y \in S${\rm ,} the distance $|x-y|$ is the square root of
an even integer. Then the subgroup of $\R^n$ generated by $S$ is
an even integral lattice.
\endproclaim

Recall that an {\it integral lattice\/} is a lattice in which
the inner product of each pair of vectors is an integer; it is
{\it even\/} if every vector has even norm.

\demo{Proof}
For all $x,y \in S$, their inner product $\langle x,y \rangle$ is
an integer, because
$$
\langle x,y \rangle = (|x-0|^2+|y-0|^2-|x-y|^2)/2 \in \Z,
$$
and the norm $|x|^2$ of any element $x \in S$ is an even integer.
It follows that the same facts hold for integer linear
combinations of elements of $S$.  The restriction on norms
implies that $S$ generates a discrete subgroup of $\R^n$, and
hence a lattice (because $S$ spans $\R^n$ by assumption).  Now
the conditions above on inner products and norms amount to the
definition of an even integral lattice.
\enddemo

Let $L$ be the subgroup of $\R^n$ generated by our periodic
packing. By Lemma~\ref{evengen}, it is an even integral lattice.
In any integral lattice, the covolume is always the square root
of an integer, since its square is the determinant of a Gram
matrix, which is an integral matrix.  Thus, $L$ has at most one
point per unit volume in $\R^n$, with equality if and only if $L$ is
unimodular. However, the periodic packing has one sphere per unit
volume in $\R^n$, because $|\Lambda|=N$.  It follows that the
periodic packing is in fact the lattice packing determined by
$L$: if any spheres from the lattice packing were missing from
the periodic packing, then by periodicity the numbers of spheres
per unit volume would be strictly smaller. Thus, our packing
comes {}from an even unimodular lattice, and that lattice must
have minimal norm $2$ in $\R^8$ and $4$ in $\R^{24}$ for the
density to be right. Such lattices are unique (see Chapters~16
and~18 of \cite{CS}). Thus, we have shown that
Conjecture~\ref{strongerconj} implies that $E_8$ and the Leech
lattice are the only periodic packings of maximal density in
$\R^8$ and $\R^{24}$, respectively.

For $\R^2$, this argument requires a slight modification, because
the isodual scaling of the hexagonal lattice is not an even
unimodular lattice.  However, the modification is not hard.  As
above, in any maximally dense periodic packing, the distances
between points must occur among the distances in the hexagonal
lattice.  Of course, we can choose any scaling we prefer;  the
most convenient is that of the $A_2$ root lattice, given by
$$
A_2 = \{(x_0,x_1,x_2) \in \Z^3 : x_0+x_1+x_2=0\},
$$
because $A_2$ is an even integral lattice.  As above, our periodic
packing must then be contained in an even integral lattice $L$.
However, the unimodularity argument no longer applies.
Fortunately, we do not need it: because $L$ is even, its minimal
norm is at least $2$, so $L$ determines a sphere packing with
spheres of the same radius as in our periodic packing.  This
sphere packing contains the original periodic packing.  If the
periodic packing did not use all these spheres, then its density
would be lower than that of $L$. Thus, it is a lattice packing,
and it is well known and easy to prove that $A_2$ is the unique
densest lattice in two dimensions.

It is worth noting that the arguments above do not require the
full strength of Conjecture~\ref{strongerconj}.  In particular,
we never make use of restrictions on the roots of the Fourier
transform.  Furthermore, in the $n=2$ case, we do not even
require as strong a condition on the function's roots: $A_2$ does
not have vectors of every even norm, but our argument allows the
function to have roots corresponding to the missing norms.

\vglue12pt\centerline{\bf Appendix A. Technicalities about density}
\vglue8pt

In this appendix, we provide precise statements and references for
what it means for a packing to have a density and whether there
is a maximally dense packing. Much of our discussion closely
follows Section~I of \cite{K}.

Let $\cal  P$ be any sphere packing in $\R^n$.  We say that
$\cal  P$ has density $\Delta$ if for all $p \in \R^n$, we have
$$
\Delta = \lim_{r \rightarrow \infty} \frac{\Vol (B(p,r) \cap
{\cal  P})}{\Vol B(p,r)},
$$
where $B(p,r)$ is the ball of radius $r$ centered at $p$, and
$B(p,r) \cap {\cal  P}$ consists of those parts of the balls
in $\cal  P$ that lie within $B(p,r)$.  It is proved in
\cite{Gr} that if this limit exists for one $p$, then it exists
for all $p$ and is equal for all $p$. We say that the packing
has uniform density if the limit exists uniformly for all $p$. In
that case, \cite{Gr} shows that for every compact set $R$ that is
the closure of its interior and every point $p$,
$$
\Delta = \lim_{r \rightarrow \infty} \frac{\Vol((rR+p) \cap
{\cal  P})}{\Vol rR},
$$
where of course $rR+p$ denotes $R$ scaled by a factor of $r$ and
translated by $p$.

Although not every packing has a density, every packing has an
upper density, defined by
$$
\Delta = \limsup_{r \rightarrow \infty} \sup_{p \in \R^n}
\frac{\Vol (B(p,r) \cap {\cal  P})}{\Vol B(p,r)}.
$$
It is proved in \cite{Gr} that the supremum of all upper densities
is achieved by a uniformly dense packing.

Periodic packings are the most convenient ones for our purposes.
Using the above results, it is easy to see that they come
arbitrarily close to the greatest possible density, as follows.
Suppose $\Delta$ is the maximum packing density in $\R^n$, and
let $\cal  P$ be a uniformly dense packing of density $\Delta$.
Let $R$ be the fundamental parallelotope of any lattice $\Lambda
\subset \R^n$.  We know that
$$
\Delta = \lim_{r \rightarrow \infty} \frac{\Vol(rR \cap {\cal
P})}{\Vol rR}.
$$
Let $\varepsilon>0$. If we choose $r$ large enough, then the
total volume of the spheres in $\cal  P$ that lie
{\it entirely\/} within $rR$ is within $\varepsilon\Vol rR$ of
$\Delta \Vol rR$, because only a negligible fraction of the
spheres can intersect the sides of $rR$.  Define a periodic
packing ${\cal  P}'$ by taking all the spheres of $\cal  P$
that lie entirely within $rR$, and also including all
translations of them by $r\Lambda$.  Then this periodic packing
has density at least $\Delta - \varepsilon$.

\vglue12pt \centerline{\bf Appendix B. Other convex bodies}
 \vglue8pt

Our methods are not limited to studying sphere packings, but
instead apply to packings with translates of any convex,
symmetrical body. In fact, it is straightforward to prove the
following generalization of our main theorem:

\specialnumber{B.1}
\proclaim{Theorem}
Let $C$ be a convex body in $\R^n${\rm ,} symmetric with
respect to the
origin.  Suppose $f : \R^n \rightarrow \R$ is an admissible
function{\rm ,} is not identically zero{\rm ,} and satisfies the following
two conditions\/{\rm :}
\begin{itemize}
\ritem{(1)} $f(x)\le 0$ for $x \not\in C${\rm ,} and
\ritem{(2)} $\widehat{f}(t) \ge 0$ for all $t$.
\end{itemize}
Then all packings with translates of $C$ have density bounded
above by
$$
\frac{\Vol(C) f(0)}{2^n{\widehat f}(0)}.
$$
\endproclaim

Unfortunately, when $C$ is not a sphere, there does not seem to
be a good analogue of the reduction to radial functions in
Theorem~\ref{main}.  That makes these cases somewhat less
convenient to deal with. There always exist functions that prove
the trivial density bound of $1$: let $\newchi_{C/2}$ be the
characteristic function of the scaled body $C/2$, and let $f =
\newchi_{C/2} * \newchi_{C/2}$.

\vglue12pt \centerline{\bf Appendix C. Numerical data}
 \vglue8pt

Table~3 compares the best packings known, the
previous best upper bound known, and our bound.  The packings and
previous bounds are from \cite{CS} (see in particular
Tables~I.1(a) and ~I.1(b) on pages~xix and~xx).
Table~4 lists the choices of forced double roots
that lead to our new bounds.

\demo{Acknowledgements}
We thank Richard Askey, Tom Brennan, Harold Diamond, Pavel
Etingof, George Gasper, David Jerison, Greg Kuperberg, Ben Logan,
L\'aszl\'o Lov\'asz, Steve Miller, Amin Shokrollahi, Neil Sloane,
Jeffrey Vaaler, David Vogan, and Michael Weinstein for helpful
discussions.\vfill
\eject
$$
\begin{tabular}{l|l|l|l}
Dimension & Best Packing Known &
Rogers' Bound & New Upper Bound\\
\hline
1&0.5     &0.5     &0.5\\
2&0.28868 &0.28868 &0.28868\\
3&0.17678 &0.1847  &0.18616\\
4&0.125   &0.13127 &0.13126\\
5&0.08839 &0.09987 &0.09975\\
6&0.07217 &0.08112 &0.08084\\
7&0.0625  &0.06981 &0.06933\\
8&0.0625  &0.06326 &0.06251\\
9&0.04419 &0.06007 &0.05900\\
10&0.03906 &0.05953 &0.05804\\
11&0.03516 &0.06136 &0.05932\\
12&0.03704 &0.06559 &0.06279\\
13&0.03516 &0.07253 &0.06870\\
14&0.03608 &0.08278 &0.07750\\
15&0.04419 &0.09735 &0.08999\\
16&0.0625  &0.11774 &0.10738\\
17&0.0625  &0.14624 &0.13150\\
18&0.07508 &0.18629 &0.16503\\
19&0.08839 &0.24308 &0.21202\\
20&0.13154 &0.32454 &0.27855\\
21&0.17678 &0.44289 &0.37389\\
22&0.33254 &0.61722 &0.51231\\
23&0.5     &0.87767 &0.71601\\
24&1.0   &1.27241 &1.01998\\
25&0.70711 &1.8798 &1.48001\\
26&0.57735 &2.8268 &2.18614\\
27&0.70711 &4.3252 &3.28537\\
28&1.0    &6.7295 &5.02059\\
29&0.70711 &10.642 &7.79782\\
30&1.0    &17.094 &12.30390\\
31&1.2095 &27.880 &19.71397\\
32&2.5658 &46.147 &32.06222\\
33&2.2220 &77.487 &52.90924\\
34&2.2220 &131.94 &88.55925\\
35&2.8284 &227.71 &150.29783\\
36&4.4394 &398.25 &258.54994\\
\end{tabular}
\begin{tabular}{c}\kern4in\end{tabular}
$$
\centerline{Table 3. Lower and upper bounds
on center density, using $m=6$.}

$$
\begin{tabular}{c|r|r|r|r|r|r|r}
Dimension & $2\pi r^2$ & $2\pi z_1^2$ & $2\pi z_2^2$ &
$2\pi z_3^2$ & $2\pi z_4^2$ & $2\pi z_5^2$ & $2\pi z_6^2$\\
\hline
 2 & 7.25520  & 21.77 & 29.02 & 50.79 & 65.34 & 90.19 & \\
 3 & 8.19385  & 22.00 & 31.63 & 49.10 & 61.80 & 80.01 & \\
 4 & 9.10543  & 22.46 & 33.49 & 48.69 & 62.03 & 79.73 & \\
 5 & 9.99512  & 23.04 & 34.87 & 48.90 & 63.04 & 80.44 & \\
 6 & 10.86682 & 23.70 & 35.98 & 49.43 & 64.07 & 81.61 & \\
 7 & 11.72351 & 24.41 & 36.94 & 50.13 & 65.04 & 82.94 & \\
 8 & 12.56674 & 25.14 & 37.74 & 50.50 & 64.08 & 79.03 & 99.37\\
 9 & 13.39945 & 25.89 & 38.57 & 51.30 & 64.99 & 80.13 & 100.44\\
10 & 14.22261 & 26.66 & 39.38 & 52.13 & 65.89 & 81.24 & 101.56\\
11 & 15.03741 & 27.44 & 40.18 & 52.97 & 66.79 & 82.34 & 102.71\\
12 & 15.84481 & 28.23 & 40.98 & 53.82 & 67.69 & 83.44 & 103.89\\
13 & 16.64565 & 29.02 & 41.77 & 54.67 & 68.61 & 84.51 & 105.08\\
14 & 17.44063 & 29.82 & 42.57 & 55.53 & 69.53 & 85.58 & 106.28\\
15 & 18.23036 & 30.61 & 43.37 & 56.39 & 70.46 & 86.65 & 107.48\\
16 & 19.01538 & 31.42 & 44.18 & 57.25 & 71.39 & 87.71 & 108.68\\
17 & 19.79617 & 32.22 & 44.99 & 58.12 & 72.34 & 88.77 & 109.88\\
18 & 20.57315 & 33.02 & 45.80 & 58.99 & 73.29 & 89.84 & 111.08\\
19 & 21.34670 & 33.82 & 46.62 & 59.86 & 74.24 & 90.90 & 112.28\\
20 & 22.11717 & 34.63 & 47.45 & 60.74 & 75.20 & 91.97 & 113.48\\
21 & 22.88488 & 35.44 & 48.27 & 61.63 & 76.17 & 93.04 & 114.68\\
22 & 23.65012 & 36.24 & 49.10 & 62.51 & 77.14 & 94.12 & 115.88\\
23 & 24.41314 & 37.05 & 49.94 & 63.41 & 78.12 & 95.20 & 117.08\\
24 & 25.17419 & 37.86 & 50.78 & 64.30 & 79.09 & 96.28 & 118.29\\
25 & 25.93349 & 38.66 & 51.62 & 65.20 & 80.08 & 97.37 & 119.50\\
26 & 26.69125 & 39.47 & 52.46 & 66.11 & 81.07 & 98.46 & 120.71\\
27 & 27.44766 & 40.28 & 53.31 & 67.02 & 82.06 & 99.55 & 121.92\\
28 & 28.20290 & 41.09 & 54.17 & 67.93 & 83.06 & 100.64 & 123.13\\
29 & 28.95711 & 41.90 & 55.02 & 68.85 & 84.06 & 101.74 & 124.34\\
30 & 29.71046 & 42.72 & 55.88 & 69.77 & 85.06 & 102.84 & 125.56\\
31 & 30.46307 & 43.53 & 56.74 & 70.69 & 86.07 & 103.94 & 126.77\\
32 & 31.21508 & 44.34 & 57.61 & 71.62 & 87.08 & 105.05 & 127.99\\
33 & 31.96659 & 45.16 & 58.47 & 72.56 & 88.09 & 106.16 & 129.21\\
34 & 32.71772 & 45.97 & 59.34 & 73.49 & 89.10 & 107.27 & 130.43\\
35 & 33.46857 & 46.79 & 60.21 & 74.43 & 90.12 & 108.38 & 131.65\\
36 & 34.21922 & 47.60 & 61.09 & 75.37 & 91.14 & 109.49 & 132.88\\
\end{tabular}
\begin{tabular}{c}\kern4in\end{tabular}
$$
\centerline{Table 4. Values of $2\pi r^2$ and $2\pi z_i^2$ for $1 \le i \le 6$.}

\AuthorRefNames [MRRW]


\begin{references}

\bibitem{AAR} \name{G.~Andrews, R.~Askey}, and \name{R.~Roy}, {\it Special
Functions\/}, Cambridge University Press, Cambridge, 1999.

\bibitem{C} \name{H.~Cohn}, New upper bounds on sphere
packings II, {\it Geom.\ Topol\/}.\  {\bf 6} (2002), 329--353,
{arXiv:math.MG/0110010}.

\bibitem{CS} \name{J.~H.~Conway} and \name{N.~J.~A.~Sloane}, {\it Sphere
Packings, Lattices and Groups\/}, third edition, Springer-Verlag, New
York, 1999.

\bibitem{D} \name{P.~Delsarte}, Bounds for unrestricted codes, by
linear programming, {\it Philips Res.\ Rep\/}.\  {\bf 27} (1972),
272--289.

\bibitem{DM} \name{H.~Dym} and \name{H.~P.~McKean}, {\it Fourier Series and
Integrals\/}, Academic Press, New York, 1972.

\bibitem{EMOT} \name{A.~Erd\'elyi, W.~Magnus, F.~Oberhettinger},
and \name{F.~G.~Tricomi}, {\it Higher Transcendental
Functions, Vol. II\/}, based, in part, on notes left by Harry
Bateman, McGraw-Hill, New York, 1953.

\bibitem{Go} \name{D.~V.~Gorbachev}, Extremal problem for entire
functions of exponential spherical type, connected with the
Levenshtein bound on the sphere packing density in $\R^n$
(Russian), {\it Izvestiya of the Tula State University Ser.\
Mathematics Mechanics Informatics\/} {\bf 6} (2000), 71--78.

\bibitem{Gr} \name{H.~Groemer}, Existenzs\"atze f\"ur Lagerungen im
Euklidischen Raum, {\it Math.\ Z\/}.\ {\bf 81}
(1963), 260--278.

\bibitem{Ha} \name{T.~Hales}, The Kepler conjecture,
1998, {arXiv:math.MG/9811078}.

\bibitem{Hs} \name{W.-Y.~Hsiang}, The optimal density of
sphere packings in dimension~$8$ and the uniqueness theorem on
final packings with optimal density, Berkeley Mathematics
Department Colloquium, March 15, 2001.

\bibitem{KL} \name{G.~A.~Kabatiansky} and \name{V.~I.~Levenshtein}, Bounds
for packings on a sphere and in space, {\it Problems of Information
Transmission\/} {\bf 14} (1978), 1--17.

\bibitem{K} \name{G.~Kuperberg}, Notions of denseness,
{\it Geom.\ Topol\/}.\  {\bf 4} (2000) 277--292,
{arXiv:\break math.MG/9908003}.

\bibitem{Leb} \name{N.~N.~Lebedev}, {\it Special Functions and Their
Applications\/}, Dover Publications, Inc., New York, 1972.

\bibitem{Lev} \name{V.~I.~Levenshtein},  Bounds for packings in
$n$-dimensional Euclidean space, {\it Soviet Math.\ Dokl\/}.\
{\bf 20} (1979), 417--421.

\bibitem{Lo} \name{B.~F.~Logan}, Extremal problems for
positive-definite bandlimited functions II: eventually negative
functions, {\it SIAM J.\ Math.\ Anal\/}.\  {\bf 14} 2 (1983), 253--257.

\bibitem{MRRW} \name{R.~J.~McEliece, E.~R.~Rodemich, H.~C.~Rumsey, Jr.,}
and
\name{L.~R.~Welch}, New upper bounds on the rate of a code via
the Delsarte-MacWilliams inequalities, {\it IEEE Trans.\
Information Theory\/} {\bf 23} (1977), 157--166.

\bibitem{OS} \name{A.~M.~Odlyzko} and \name{N.~J.~A.~Sloane}, New bounds
on the number of unit spheres that can touch a unit sphere in
$n$~dimensions, {\it J.\ Combin.\ Theory  A\/} {\bf 26}
(1979), 210--214.

\bibitem{Ro} \name{C.~A.~Rogers}, The packing of equal spheres,
{\it Proc.\  London Math.\  Soc\/}.\  {\bf 8} (1958),
609--620.

\bibitem{Ru} \name{W.~Rudin}, {\it Real and Complex Analysis\/},
third edition, McGraw-Hill, New York, 1987.

\bibitem{S} \name{C.~L.~Siegel}, \"Uber Gitterpunkte in convexen
K\"orpern and ein damit zusammenh\"angendes Extremalproblem,
{\it Acta Math\/}.\  {\bf 65} (1935), 307--323.

\bibitem{SW} \name{E.~M.~Stein} and \name{G.~L.~Weiss},
{\it Introduction to Fourier Analysis on Euclidean Spaces},
{\it Princeton Math.\  Series\/} {\bf 32}, Princeton Univ.\  Press,
Princeton, NJ, 1971.

\bibitem{Y} \name{V.~A.~Yudin}, The Laplace operator and  packings
of balls  in a Euclidean space, {\it Differential Equations\/} {\bf 31}
 (1995), 837--840.
\end{references}
\end{document}